\documentclass[]{interact}

\usepackage{epstopdf}
\usepackage[caption=false]{subfig}

\usepackage[numbers,sort&compress]{natbib}
\bibpunct[, ]{[}{]}{,}{n}{,}{,}
\makeatletter
\def\NAT@def@citea{\def\@citea{\NAT@separator}}
\makeatother
\usepackage{anyfontsize}
\theoremstyle{plain}
\newtheorem{theorem}{Theorem}[section]
\newtheorem{lemma}[theorem]{Lemma}
\newtheorem{corollary}[theorem]{Corollary}
\newtheorem{proposition}[theorem]{Proposition}

\theoremstyle{definition}
\newtheorem{definition}[theorem]{Definition}
\newtheorem{remark}[theorem]{Remark}
\newtheorem{example}[theorem]{Example}

\usepackage{xcolor}

\theoremstyle{remark}

\begin{document}
	

	\title{A formula of $A$-spectral radius for $A^{\frac{1}{2}}$-adjoint operators on semi-Hilbertian spaces}
  \author{\name {Arup Majumdar \thanks{Arup Majumdar (corresponding author). Email address: arupmajumdar93@gmail.com}and P. Sam Johnson \thanks{ P. Sam Johnson. Email address: sam@nitk.edu.in}} \affil{Department of Mathematical and Computational Sciences, \\
			National Institute of Technology Karnataka, Surathkal, Mangaluru 575025, India.}}
  
	\maketitle

	\begin{abstract}
		In this paper, we prove the relation $\frac{r_{A}(T) + r_{A}(T^{\diamond}) + |r_{A}(T^{\diamond}) - r_{A}(T)|}{2} = \sup \{ |\lambda|: \lambda \in \sigma_{A}(T)\}$, where $A$ is a positive semidefinite operator (not necessarily to have a closed range) and $r_{A}(T)$ is the $A$-spectral radius of $T$ in $B_{A^{\frac{1}{2}}}(H)$. Also we prove that $\sup \{ |\lambda|: \lambda \in \sigma_{A}(T)\} = r_{A}(T), \text{ when }  T \in B_{A^{\frac{1}{2}}}(H)  \text { commutes with } A$. By introducing $A$-Harte spectrum $\sigma_{A_{h}}(\bold{T})$ of a $d$-tuple operator $\bold{T}= (T_{1},\dots,T_{d}) \in (B_{A^{\frac{1}{2}}}(H))^{d}$, we prove that  $r_{A_{h}}(\bold{T}) \leq \sup \{\|\lambda\|_{2}: \lambda \in \sigma_{A_{h}}(\bold{T})\}$, where $r_{A_{h}}(\bold{T})$ is the $A$-Harte spectral radius of $\bold{T}$.
	\end{abstract}
	
      \begin{amscode}47A10; 46C05; 47A30 ; 47A80.\end{amscode}
    \section{Introduction}
    In quantum mechanics, the Hamiltonian operator is self-adjoint with respect to the chosen inner product. Therefore, a non-Hermitian operator, $T$, has to be Hermitian with respect to a different inner product to serve as the Hamiltonian for a unitary quantum system. The most general inner product that makes $T$ as Hermitian has the form $\langle x, y \rangle_{A} = \langle Ax, y \rangle$, where $\langle.,. \rangle$ is the defining inner product on the original Hilbert space, and $A$ is a positive-definite operator satisfying the pseudo-Hermiticity condition $T^{*} = ATA^{-1}$ {\cite{{MR1144600, MR2252709}}. The operator $A$ is called the metric operator, which is the basic ingredient of quantum theory based on the Hamiltonian operator $T$. The larger class of quasi-Hermitian operators gives useful flexibility in the mathematical description of physical observables in quantum mechanics. In the mathematical approach, $A$ is considered as a positive semidefinite operator and new semi-inner product space is considered with the semi-inner product, $\langle x, y \rangle_{A} = \langle Ax, y \rangle$,  called the semi-Hilbertian space. In the literature, many papers study $A$-bounded operators in semi-Hilbertian spaces.
    
    This paper is devoted to studying the $A$-approximate spectrum of $A$-bounded operators acting on a complex semi-Hilbertian space $H$. In 2008, Arias et al. investigated partial isometries and metric properties of projections in semi-Hilbertian spaces \cite{MR2442900, MR2388631}. After that, many mathematicians have been working on semi-Hilbertian spaces to investigate whether the basic properties of bounded operators on Hilbert spaces will work or not for operators in semi-Hilbertian spaces. Kais Feki introduced the notion of the $A$-spectral radius of the class of  operators in $B_{A^{\frac{1}{2}}}(H)$ acting  on a complex Hilbert space $H$ \cite{MR4133606}. Hamadi et al. have developed the spectral analysis of $A$-bounded operator $T \in B_{A^{\frac{1}{2}}}(H)$  and established some spectral results of those operators $T \in B_{A^{\frac{1}{2}}}(H)$ with the condition that the range of $A$ is closed \cite{MR4362420, hb}. 

    The paper is organized as follows. We start with some notations in section 2. In section 3,  without the assumption of closed range of $A$, we prove the relation $\frac{r_{A}(T) + r_{A}(T^{\diamond}) + |r_{A}(T^{\diamond}) - r_{A}(T)|}{2} = \sup \{ |\lambda|: \lambda \in \sigma_{A}(T)\}$, when $T \in B_{A^{\frac{1}{2}}}(H)$. Also we prove that $\sup \{ |\lambda|: \lambda \in \sigma_{A}(T)\} = r_{A}(T)$ when $T \in B_{A^{\frac{1}{2}}}(H) \text { commutes with } A$. In the final section, we introduce $A$-Harte spectrum $\sigma_{A_{h}}(\bold{T})$ of a $d$-tuple operator $\bold{T}= (T_{1},\dots,T_{d}) \in (B_{A^{\frac{1}{2}}}(H))^{d}$ and prove that  $r_{A_{h}}(\bold{T}) \leq \sup \{\|\lambda\|_{2}: \lambda \in \sigma_{A_{h}}(\bold{T})\}$, where $r_{A_{h}}(\bold{T})$ is the $A$-Harte spectral radius of $\bold{T}$.
   
	\section{Preliminaries}
 Throughout the paper, $H$ denotes a Hilbert space over a complex field. The algebra of all linear bounded operators on $H$ is denoted by $B(H)$, whereas $B(H)^{+}$ is the cone of positive semidefinite operators, i.e., $B(H)^{+} = \{ T\in B(H) : \langle Tx, x\rangle \geq 0\, \text { for all } x\in H\}$. The operator $A$ is a non-zero positive semidefinite operator. The sesquilinear form defines from $H \times H$ to $\mathbb C$  by ${\langle \xi , \eta \rangle}_{A} = \langle A\xi , \eta \rangle$. Moreover, ${\| .\|}_{A}$ indicates the seminorm induced by ${\langle . ,. \rangle}_{A}$, that is, ${\| \xi \|}_{A} = {\langle \xi , \xi \rangle}_{A}^{\frac{1}{2}}$. It is obvious that $\| \xi\|_{A} = 0$ if and only if $\xi \in N(A)$, where $N(A)$ is the null space of $A$.  The subspace ${\mathfrak{S}}^{\perp_{A}} = \{\xi : \langle A\xi, \eta \rangle = 0 \text{ for all } \eta \in \mathfrak{S}\}$ is called the $A$-orthogonal companion of $\mathfrak{S}$. The range of an operator $T$ is denoted by $R(T)$ and its closure is denoted by  $\overline{R(T)}$.
 \begin{definition}\cite{MR2388631} An operator  $T \in B(H)$ is called as $A$-bounded operator if there exists a constant $c > 0$ such that $\|T \xi\|_{A} \leq c  \|\xi\|_{A}$,    for all $\xi \in \overline{ R(A)}$. Moreover,
 \begin{equation*}
     {\|T\|}_{A} = \sup_{\xi \in \overline {R(A)}\setminus \{0\}}\frac{{\|T \xi\|}_{A}}{{\| \xi\|}_{A}} < \infty.
 \end{equation*}
 Equivalently, it can be observed that 
 \begin{align*}
 {\|T\|}_{A} = \sup \{ | {\langle T \xi, \eta \rangle}_{A} | : \xi, \eta \in H , {\| \xi\|}_{A} \leq 1, {\| \eta\|}_{A} \leq 1\}.
 \end{align*}
 We denote
 \begin{align*}
 B^{A}(H)= \{T \in B(H) : {\|T\|}_{A} < \infty\}.
 \end{align*}
 
 \end{definition}
 \begin{theorem} \cite{MR0203464}\label{thm 1}
 Let $E, F \in B(H)$. The following conditions are equivalent:
 \begin{enumerate}
     \item $R(F) \subset R(E)$.
     \item There exists a positive number $\lambda $ such that $FF^{*} \leq \lambda EE^{*}$.
     \item There exists $C\in B(H)$ such that $EC = F$.
 \end{enumerate}
 If one of these conditions holds, then there exists a unique operator $D \in B(H)$ such that $ED = F$ with $R(D) \subset \overline {R(E^{*})}$. Furthermore, $N(D) = N(F)$. The unique operator, $D$, is called the reduced solution of $EX = F$.
 \end{theorem}
 \begin{definition}\cite{MR2442900}
 A bounded linear operator $W$ is called an $A$-adjoint of $T\in B(H)$ if 
 \begin{align*}
 {\langle T\xi , \eta \rangle}_{A} = {\langle \xi, W\eta \rangle}_{A} ,  \text{ for all }   \xi, \eta \in H.
 \end{align*}
 \end{definition}
 From Theorem \ref{thm 1}, the existence of an $A$-adjoint operator of $T\in B(H)$ is guaranteed if and only if $R(T^{*}A) \subset R(A)$. By $B_{A}(H)$, we denote the subalgebra of $B(H)$ consisting of such operators which have $A$-adjoint operators, that is,
 \begin{align*}
 B_{A}(H) = \{T \in B(H) : R(T^{*}A) \subset R(A)\}.
 \end{align*}
 Similarly, we can define ${B_{A^{\frac{1}{2}}}}(H) = \{ T \in B(H) : R(T^{*}A^{\frac{1}{2}}) \subset R(A^{\frac{1}{2}})\}$. Again, by Theorem \ref{thm 1}, it is obvious that 
 \begin{align*}
 B_{A^{\frac{1}{2}}}(H) = \{ T \in B(H) : \exists ~c>0, ~{\|T\xi\|}_{A} \leq c {\|\xi\|}_{A},  \text { for all } \xi \in H\}.
 \end{align*}
We have the  following inclusion relations \cite{MR2590353} : $$B_{A}(H) \subset B_{A^{\frac{1}{2}}}(H) \subset B^{A}(H) \subset B(H).$$
 
If $T\in B_{A}(H)$, then there exists a unique A-adjoint, say $T^{\sharp}$, such that $T^{*}A = A T^{\sharp}$ with $R(T^{\sharp}) \subset \overline {R(A)}$. The reduced solution $T^{\sharp}$ of $T^{*}A= AX$ gives that $N(T^{\sharp}) = N(T^{*}A)$ and $T^{\sharp} = A^{\dagger}T^{*}A$, where $A^{\dagger}$ denotes the Moore-Penrose of $A$ in the domain $ D(A^{\dagger}) =R(A) \oplus R(A)^{\perp}$ with 
\begin{align*}
A^{\dagger} y = 
			\begin{cases}
				(A_{N(A)^{\perp}})^{-1}y  & \text{if} ~  y\in R(A)\\
				0    & \text{if}  ~ y\in R(A)^{\perp}.
			\end{cases}
\end{align*}
Similarly, if $T \in B_{A^{\frac{1}{2}}}(H)$, there exists a unique reduced solution $T^{\diamond}$ such that $T^{*}A^{\frac{1}{2}} = A^{\frac{1}{2}} T^{\diamond}$ with $R(T^{\diamond}) \subset \overline{R(A^{\frac{1}{2}})}$.
\begin{proposition}\cite{MR2442900}\label{ pro 2.5}
Let $T\in B(H)$. Then the following conditions are equivalent:
\begin{enumerate}
    \item $T \in B^{A}(H)$.
    \item $A^{\frac{1}{2}}T(A^{\frac{1}{2}})^{\dagger}$ is a bounded linear operator.
    \item $R(A^{\frac{1}{2}}T^{*}A^{\frac{1}{2}}) \subset R(A)$.
\end{enumerate}
Moreover, one of the above conditions guarantees that ${\|T\|}_{A} = \|A^{\frac{1}{2}}T(A^{\frac{1}{2}})^{\dagger}\| = \|(A^{\frac{1}{2}})^{\dagger} T^{*}A^{\frac{1}{2}}\|$.
\end{proposition}
Some properties of $T^{\sharp}$ have been studied in the papers \cite{MR2442900, MR2388631}. For the sake of completeness, we mention some of them. From now on, we write $P$ instead of $P_{\overline{R(A)}}$. Here, $P_{\overline{R(A)}}$ denotes the orthogonal projection onto $\overline{R(A)}$.
\begin{proposition}\cite{MR2442900, MR2388631}
    Let $T \in B_{A}(H)$. Then the following statements are true:
    \begin{enumerate}
        \item $(A^{a})^{\sharp} = A^{a}$ for every $a > 0$.
        \item If $AT = TA$, then $T^{\sharp} = PT^{*}$.
        \item If $AT = T^{*}A$, then $(A^{\frac{1}{2}})^{\dagger}T^{*}A^{\frac{1}{2}}$ is selfadjoint.
        \item If $W \in B_{A}(H)$, then $TW \in B_{A}(H)$ and $(TW)^{\sharp}= W^{\sharp}T^{\sharp}$.
        \item $T^{\sharp} \in B_{A}(H)$, $(T^{\sharp})^{\sharp} = PTP$ and $((T^{\sharp})^{\sharp})^{\sharp} = T^{\sharp}$.
        \item ${\|T\|}_{A} = {\|T^{\sharp}\|}_{A} = {{\|T^{\sharp}T\|}_{A}}^{\frac{1}{2}}$.
        \item ${\|W\|}_{A} = {\|T^{\sharp}\|}_{A}$ for every $A$-adjoint $W \in B(H)$ of $T$.
        \item $\|T^{\sharp}\| \leq \|W\|$ for every $W\in B(H)$ which is an $A$-adjoint of $T$. Nevertheless, $T^{\sharp}$ is not generally unique $A$-adjoint of $T$ that has the minimal norm.
\end{enumerate}
\end{proposition}
\begin{definition}\cite{MR2388631} An operator      $T \in B(H)$ is called an $A$-isometry if ${\|T \xi\|}_{A} = {\|\xi\|}_{A}  , \text{ for all } \xi \in H$.
\end{definition}
\noindent It is shown in \cite{MR2388631} that  for  $T \in B_{A}(H)$,  $T$ is an $A$-isometry if and only if $T^{\sharp}T = P$.

\begin{definition}\cite{MR2388631}
    Let $U\in B_{A}(H)$. Then $U$ is called $A$-unitary if $U$ and $U^{\sharp}$ both are $A$-ismotries.
\end{definition}
\begin{proposition}\cite{MR2388631}
    Let $U, V \in B_{A}(H)$ be $A$-unitary operators. Then
    \begin{enumerate}
        \item $U^{\sharp}U = (U^{\sharp})^{\sharp} U^{\sharp} = P$.
        \item ${\|U\|}_{A}= 1$.
        \item $U^{\sharp}$ is $A$-unitary.
        \item $UV$ is $A$-unitary.
        \item ${\|UTV^{\sharp}\|}_{A} = {\|T\|}_{A}$, for every $T \in B_{A^{\frac{1}{2}}}(H)$.
    \end{enumerate}
\end{proposition}
\begin{definition}\cite{MR4133606}
Let $T \in B(H)$. The $A$-numerical range is defined by 
\begin{align*}
W_{A}(T) = \{\langle Tx, x\rangle_{A} : x\in H, \|x\|_{A} = 1\}.
\end{align*}
Moreover, $w_{A}(T) = \sup \{|\langle Tx, x\rangle|_{A} : x\in H, \|x\|_{A}= 1\}$ is called the $A$-numerical radius of $T$.
\end{definition}
\begin{theorem}\cite{MR3834203}\label{thm 2.14}
Let $T \in B(H)$. Then $W_{A}(T)$ is a convex subset of $\mathbb C$.
\end{theorem}

 \begin{definition}\cite{MR4362420}
  A non-zero operator $T\in B_{A^{\frac{1}{2}}}(H)$ is said to be $A$-invertible in $ B_{A^{\frac{1}{2}}}(H)$ if there exists a non-zero operator $S\in B_{A^{\frac{1}{2}}}(H)$ such that $ATS = AST =A$. Here $S$ is called an $A$-inverse in $ B_{A^{\frac{1}{2}}}(H) $.  
  
  In a similar way, we say that a non-zero operator $T \in  B_{A}(H)$ is said to be $A$-invertible in $ B_{A}(H)$ if there exists a non-zero operator $S \in  B_{A}(H)$ such that $ATS = AST =A$. Here $S$ is called an $A$-inverse in $B_{A}(H)$.
 \end{definition}
 \begin{definition}\cite{MR4362420}
 Let $T \in B_{A^{\frac{1}{2}}}(H)$.
 \begin{enumerate}
     \item The set $\rho_{A}(T) = \{ \lambda \in \mathbb {C}: (\lambda - T)$ is $A$-invertible in $B_{A^{\frac{1}{2}}}(H)\}$ is called the $A$-resolvent set of $T$.
     \item The set $\sigma_{A}(T) = \mathbb{C} \setminus \rho_{A}(T)$ is called the $A$-spectrum of $T$.
 \end{enumerate}
 \end{definition}
  \begin{remark}\cite{MR4362420}
 Let $T \in B_{A^{\frac{1}{2}}}(H)$ be an $A$-invertible operator in $B_{A^{\frac{1}{2}}}(H)$ with an $A$-inverse $S \in B_{A^{\frac{1}{2}}}(H)$. The following statements are equivalent:
 \begin{enumerate}
   \item $ATS = AST =A$.
   \item $PTS = PST =P$.
   \item $A^{\frac{1}{2}}TS = A^{\frac{1}{2}}ST = A^{\frac{1}{2}}$.
 \end{enumerate}
 \end{remark}
 \begin{remark}\label{remark 1}
 Let $T, S$ be in $B_{A}(H)$. The operator $T$ is an $A$-invertible operator in $B_{A}(H)$ with an $A$-inverse $S \text { in } B_{A}(H)$ if and only if $T^{\sharp}$ is also $A$-invertible in $B_{A}(H)$ with an $A$-inverse $S^{\sharp} \text { in } B_{A}(H)$.
 \end{remark}
\begin{definition}\cite{MR4133606}
Let $T \in B_{A^{\frac{1}{2}}}(H)$. The $A$-spectral radius of $T$ is defined as
\begin{align*}
r_{A}(T) = \lim_{n \to \infty} (\|T^{n}\|_{A})^{\frac{1}{n}}.
\end{align*}
\end{definition}
\begin{theorem}\cite{MR4133606}\label{thm 2.15}
Let $T \in B_{A^{\frac{1}{2}}}(H)$. Then $r_{A}(T) \leq w_{A}(T) \leq \|T\|_{A}$.
\end{theorem}
 \begin{theorem}\cite{hb}\label{thm 2}
 If $T \in B_{A}(H)$  is $A$-invertible in $B_{A^{\frac{1}{2}}}(H)$, then each $A$-inverse of $T$ in $B_{A^{\frac{1}{2}}}(H)$ is in $B_{A}(H)$.
 \end{theorem}
\begin{definition}\cite{MR2590353}
Let $T\in B(H)$. The $A$-reduced minimum modulus of $T$ is given by 
    \begin{align}
    {\gamma_{A}}(T) = \inf \{\|T\xi\|_{A}: \xi \in {N(A^{\frac{1}{2}}T)^{\perp_{A}} , \|\xi\|_{A}=1}\}.
    \end{align}
Note that if $T \in B_{A}(H)$, then ${\gamma_{A}}(T) =  \inf \{\|T\xi\|_{A}: \xi \in \overline{R(T^{\sharp}T)} , \|\xi\|_{A}=1\}$.
\end{definition}
\begin{proposition}\cite{MR2590353}
Let $T \in B_{A}(H)$ and $C$ be a solution of the equation $A^{\frac{1}{2}}X = T^{*}A^{\frac{1}{2}}$. If $A^{\frac{1}{2}} \overline{R(T^{\sharp}T)} \subset \overline {R(C)}$, then $\gamma_{A}(T) = \gamma(C)$.
\end{proposition}
\begin{proposition}\cite{MR2590353}\label{pro 0.1}
Let $T \in B_{A}(H)$. Then
\begin{enumerate}
 \item $\gamma_{A}(T) = \gamma(T^{\diamond})$.
 \item $\gamma_{A}(T) = \gamma_{A}(T^{\sharp})$.
\end{enumerate}
\end{proposition}
 Let us now consider the Hilbert space $\mathbf{R}(A^{\frac{1}{2}}) = (R(A^{\frac{1}{2}}), (.,.))$  with the inner product $(A^{\frac{1}{2}} \xi, A^{\frac{1}{2}} \eta) = \langle P\xi, P\eta\rangle,  \text {for all } \xi, \eta \in H$. Then $\|A^{\frac{1}{2}} \xi\|_{\mathbf{R}(A^{\frac{1}{2}})} = \|P\xi\|, \text { for all } \xi \in H$. We consider an operator $W_{A} : H \mapsto \mathbf{R}(A^{\frac{1}{2}})$ defined by $W_{A}(\xi) = A\xi, \text { for all } \xi\in H$ and $Z_{A}: H \mapsto \mathbf{R}(A^{\frac{1}{2}})$ defined by $Z_{A}(\xi) = A^{\frac{1}{2}}(\xi), \text{ for all } \xi\in H$.
 \begin{proposition}\cite{MR2590353}
  The following assertions hold:
  \begin{enumerate}
   \item $Z_{A} \in B(H,\mathbf{R}(A^{\frac{1}{2}}))$ and $Z_{A}$ is onto.
   \item $Z_{A}^{*} \in B(\mathbf{R}(A^{\frac{1}{2}}), H)$, $Z_{A}^{*}(A^{\frac{1}{2}}(\xi))= P(\xi)$.
   \item $Z_{A}^{*}Z_{A}= P$ and $Z_{A}Z_{A}^{*}= I_{\mathbf{R}(A^{\frac{1}{2}})}$.
   \item $W_{A} \in B(H, \mathbf{R}(A^{\frac{1}{2}}))$ and $R(W_{A}) = R(A)$ is dense in $\mathbf{R}(A^{\frac{1}{2}})$.
   \item $W_{A}^{*}(A^{\frac{1}{2}}(\xi)) = A^{\frac{1}{2}}(\xi)$, and $R(W_{A}^{*}) = R(A^{\frac{1}{2}})$.
  \end{enumerate}
 \end{proposition}
 
 The following result describes the relationship among $A$-bounded operators in a semi-Hilbertian space $H$ with the operators of $B(\mathbf{R}(A^{\frac{1}{2}}))$. 
 
\begin{proposition}\cite{MR2590353}
Let $T \in B(H)$. Then there exists $\tilde{T} \in B(\mathbf{R}(A^{\frac{1}{2}}) )$ such that $\tilde{T}W_{A} = W_{A}T$ if and only if $T\in B_{A^{\frac{1}{2}}}(H)$. In this case, $\tilde{T}$ is unique.
\end{proposition}
 A new mapping $\alpha: B_{A^{\frac{1}{2}}}(H) \mapsto B(\mathbf{R}(A^{\frac{1}{2}}) )$ is defined as $\alpha(T) = \tilde{T}$. The next proposition divulges some properties of $\alpha$.
 \begin{proposition}\cite{MR2590353}
 The following assertions hold:
 \begin{enumerate}
     \item $\alpha$ is the homomorphism and $\alpha(T)= \overline {W_{A}T W_{A}^{\sharp}}$, where $W_{A}^{\sharp} = W_{A}^{\dagger}$.
     \item $\|\alpha(T)\|_{\mathbf{R}(A^{\frac{1}{2}})} = \|T\|_{A}$.
 \end{enumerate}
 \end{proposition}
\section{Some characterizations of the $A$-spectral radius of $A^{\frac{1}{2}}$-adjoint operators}
Hamadi Baklouti et al. raised an open question in \cite{MR4362420} whether $\sup\{ |\lambda|: \lambda \in \sigma_{A}(T)\} = r_{A}(T)$ holds good for $T \in B_{A^{\frac{1}{2}}}(H)$. The proof of Theorem $5.1$ in \cite{MR4362420} requires $R(A)$ to be closed, though it is not mentioned in the statement. In the same proof, one can see that  inequality  $\|((A^{\frac{1}{2}})^{\dagger} (\lambda I -T)^{*} A^{\frac{1}{2}})x\| \geq (\lambda -\|T\|_{A})\|Px\| ~ ( \text{ for all } x\in H)$ would be proper instead of the inequality $\|((A^{\frac{1}{2}})^{\dagger} (\lambda I -T)^{*} A^{\frac{1}{2}})x\| \geq (\lambda -\|T\|_{A})\|x\| ~( \text{ for all } x\in H)$. Furthermore, Kais Feki has proved $\sup \{ |\lambda|: \lambda \in \sigma_{A}(T)\} \geq r_{A}(T), \text{ when } T \in B_{A^{\frac{1}{2}}}(H)$, in Theorem 3.1 \cite{MR4402703}. However, the equality part of Theorem 3.1 \cite{MR4402703} is not exact because the author has considered the result of Theorem 5.1 \cite{MR4362420} to show $\sup \{ |\lambda|: \lambda \in \sigma_{A}(T)\} \leq r_{A}(T)$ without assuming $R(A)$ to be closed. Therefore, the problem $\sup \{ |\lambda|: \lambda \in \sigma_{A}(T)\} = r_{A}(T), \text{ when } T \in B_{A^{\frac{1}{2}}}(H)$ is still open.  

In this section, we establish the relation $\frac{r_{A}(T) + r_{A}(T^{\diamond}) + |r_{A}(T^{\diamond}) - r_{A}(T)|}{2} = \sup \{ |\lambda|: \lambda \in \sigma_{A}(T)\}$. Also we prove that $\sup \{ |\lambda|: \lambda \in \sigma_{A}(T)\} = r_{A}(T), \text{ when }  T \in B_{A^{\frac{1}{2}}}(H)  \text { commutes with } A$.

 We have corrected the statement of Theorem 5.1 \cite{MR4362420} in Theorem \ref{thm 3} by adding the assumption of $R(A)$ to be closed. 
\begin{theorem}\label{thm 3}
     Let $T \in B_{A^{\frac{1}{2}}}(H)$ with $R(A)$ is closed. Then for all $|\lambda | > \|T\|_{A}$, $\lambda \in \rho_{A}(T)$.
 \end{theorem}
 \begin{proof}
 It is proved in Proposition \ref{ pro 2.5} that $\|T\|_{A} = \| A^{\frac{1}{2}} T  (A^{\frac{1}{2}})^{\dagger}\| = \|(A^{\frac{1}{2}})^{\dagger} T^{*} A^{\frac{1}{2}}\|$. Moreover, we can easily show that 
 \begin{equation}\label{ineq 1}
     \|(A^{\frac{1}{2}})^{\dagger} (\lambda -T)^{*} A^{\frac{1}{2}}x\| \geq (|\lambda| - \|T\|_{A}) \|Px\|, \  \text{ for all }\ x\in H.
 \end{equation}
Since $R(A)$ is closed. So, $T \in B_{A}(H) $. The above inequality (\ref{ineq 1}) implies,

\begin{align*}
    \begin{split}
         & ~~~~{A^{\frac{1}{2}} (\lambda -T) (A^{\frac{1}{2}})^{\dagger}} (A^{\frac{1}{2}})^{\dagger} (\lambda -T)^{*} A^{\frac{1}{2}} \geq (|\lambda| - \|T\|_{A})^{2} P \\
         & =A(\lambda -T)(A^{\frac{1}{2}})^{\dagger}(A^{\frac{1}{2}})^{\dagger}(\lambda -T)^{*}A \geq (|\lambda| - \|T\|_{A})^{2} A.
        \end{split}
\end{align*}

By Theorem \ref{thm 1}, there exists a $W\in B(H)$ such that $A(\lambda -T)(A^{\frac{1}{2}})^{\dagger}W A^{\frac{1}{2}} = A$. Thus, $R(A) \subset R(A(\lambda -T))$. The reverse inclusion is obvious, which gives that $R(A)= R(A(\lambda -T))$.
It is straightforward to prove that 
\begin{align*}
(|\lambda| - \|T\|_{A}) \|x\|_{A} \leq \|(T- \lambda) x\|_{A} \leq (\|T\|_{A} + |\lambda|)\|x\|_{A}, \  \text{ for all }\ x\in H.
\end{align*}
From Theorem 4.2 in \cite{MR4362420} and Theorem \ref{thm 2}, we conclude that the operator $(\lambda -T)$ is invertible in $B_{A}(H)$. Thus, $\lambda \in \rho_{A}(T)$.
\end{proof}
\begin{remark} Let $T \in  B_{A^{\frac{1}{2}}}(H)$ with $R(A)$ is closed. Then Theorem \ref{thm 3} and the first part of Proposition 5.10 \cite{MR4362420} say that $\sup\{|\lambda| : \lambda \in \sigma_{A}(T) \} \leq r_{A}(T)$.

\end{remark}
The following Theorem is proved by Kais Feki \cite{MR4402703} but we present a different approach to prove that $\sup \{ |\lambda|: \lambda \in \sigma_{A}(T)\} \geq r_{A}(T) , \text{ when } T \in B_{A^{\frac{1}{2}}}(H)$.
\begin{theorem}\label{thm 5.1}
Let $T \in B_{A^{\frac{1}{2}}}(H)$. Then $\sup \{ |\lambda|: \lambda \in \sigma_{A}(T)\} \geq r_{A}(T)$.
\end{theorem}
\begin{proof}
We define a new operator $T_{a}: \overline{R(A)} \mapsto \overline{R(A)}$ by $T_{a}(x) = (A^{\frac{1}{2}})^{\dagger} T^{*} A^{\frac{1}{2}} x, \text { for all } x \in \overline{R(A)}$. Now consider $\lambda \in \rho_{A}(T)$, then there exists $S_{\lambda} \in B_{A^{\frac{1}{2}}}(H)$ such that $P(\lambda -T) S_{\lambda} = P S_{\lambda} (\lambda -T) = P$. So,
\begin{align*}
 & (S_{\lambda})^{*} (\overline{\lambda} - T^{*})P = (\overline{\lambda} - T^{*}) (S_{\lambda})^{*}P = P \\
&\implies (A^{\frac{1}{2}})^{\dagger} (S_{\lambda})^{*} (\overline{\lambda} - T^{*})A^{\frac{1}{2}} = (A^{\frac{1}{2}})^{\dagger} (\overline{\lambda} - T^{*}) (S_{\lambda})^{*} A^{\frac{1}{2}} = I_{\overline{R(A)}}\\
& = (A^{\frac{1}{2}})^{\dagger} (S_{\lambda})^{*} A^{\frac{1}{2}} (A^{\frac{1}{2}})^{\dagger} (\overline{\lambda} - T^{*})A^{\frac{1}{2}} = (A^{\frac{1}{2}})^{\dagger} (\overline{\lambda} - T^{*}) A^{\frac{1}{2}} (A^{\frac{1}{2}})^{\dagger} (S_{\lambda})^{*} A^{\frac{1}{2}} = I_{\overline{R(A)}}\\
& = (S_{\lambda})_{a} (\overline{\lambda} -T_{a}) = (\overline{\lambda} - T_{a}) (S_{\lambda})_{a} = I_{\overline{R(A)}}.
\end{align*}
Here, $(S_{\lambda})_{a} = (A^{\frac{1}{2}})^{\dagger} (S_{\lambda})^{*} A^{\frac{1}{2}}$ on domain $\overline{R(A)}$. Thus, $\overline{\lambda} \in \rho(T_{a})$. Hence 
\begin{align}\label{ineq 17}
\lim_{n \to \infty}\|(T_{a})^{n}\|^{\frac{1}{n}} = \sup \{|\lambda| : \lambda \in \sigma(T_{a})\} \leq \sup \{|\lambda| : \lambda \in \sigma_{A}(T)\}.
\end{align}
We know that  $\|T_{a}\| = \|(A^{\frac{1}{2}})^{\dagger} T^{*}A^{\frac{1}{2}}\| = \|T\|_{A}$, by Proposition \ref{ pro 2.5}. So, $\|(T_{a})^{n}\| = \| (A^{\frac{1}{2}})^{\dagger} (T^{n})^{*} A^{\frac{1}{2}}\|= \| T^{n}\|_{A}, \text { for all }  n\in \mathbb N$. Therefore, the inequality (\ref{ineq 17}) gives that $r(T_{a}) = r_{A}(T) \leq \sup \{|\lambda| : \lambda \in \sigma_{A}(T)\}$, where $r(T_{a})$ is the spectral radius of $T_{a}$.
\end{proof}
\begin{remark}
    Let $T \in B_{A^{\frac{1}{2}}}(H)$ with $R(A)$ is closed. Then $\sup \{ |\lambda|: \lambda \in \sigma_{A}(T)\} = r_{A}(T)$.
\end{remark}
   
For $T\in B(H)$, the effective action $T_{eff}$ of $T$ with respect to $A$ is defined on  $\overline{R(A)}$ by $T_{eff}(x) = PTP(x), \text { for all } x\in \overline{R(A)}$.
\begin{lemma}\label{lemma 5} \cite{hb}
Let $T \in B_{A^{\frac{1}{2}}}(H)$. If $T$ is $A$-invertible in $B_{A^{\frac{1}{2}}}(H)$, then $T_{eff}$ is invertible in $B(\overline{R(A)})$. 
\end{lemma}
\begin{lemma}\label{lemma 0.3.4}
Let $T \in B_{A^{\frac{1}{2}}}(H)$. Then $r(T_{eff}) = r_{A}(T^{\diamond}) \leq \sup \{|\lambda|: \lambda \in \sigma_{A}(T)\}$.
\end{lemma}
\begin{proof}
By Lemma \ref{lemma 5}, it is obvious that $r(T_{eff}) = \sup \{|\lambda| : \lambda \in \sigma(T_{eff})\} \leq \sup \{ |\lambda|: \lambda \in \sigma_{A}(T)\}$.
 
 \begin{align*}
\|T_{eff}\| &= \|PTP\| \\
&= \|(A^{\frac{1}{2}})^{\dagger} A^{\frac{1}{2}} T (A^{\frac{1}{2}})^{\dagger} A^{\frac{1}{2}}\| \\
&= \|(A^{\frac{1}{2}})^{\dagger} \overline{A^{\frac{1}{2}} T (A^{\frac{1}{2}})^{\dagger}} A^{\frac{1}{2}}\| \\
&= \|(A^{\frac{1}{2}})^{\dagger} (A^{\frac{1}{2}} T (A^{\frac{1}{2}})^{\dagger})^{**} A^{\frac{1}{2}}\| \\
&= \|(A^{\frac{1}{2}} T (A^{\frac{1}{2}})^{\dagger})^{*}\|_{A} \\
&= \|(A^{\frac{1}{2}})^{\dagger} T^{*} A^{\frac{1}{2}}\|_{A}.
\end{align*}
Similarly, we can prove that $\|(T_{eff})^{n}\| = \|(A^{\frac{1}{2}})^{\dagger} (T^{n})^{*} A^{\frac{1}{2}}\|_{A}, \text{ for all } n \in \mathbb N$. Now, $\|(A^{\frac{1}{2}})^{\dagger} (T^{n})^{*} A^{\frac{1}{2}} x\|_{A} = \|(T^{n})^{*} A^{\frac{1}{2}} x\| = \|(T^{n})^{\diamond}x\|_{A}, \text{ for all } n\in \mathbb N$ and $x\in \overline{R(A)}$. Thus, $\|(T_{eff})^{n}\| = \| (T^{\diamond})^{n}\|_{A}, \text{ for all }n \in \mathbb N$. Therefore,  $r(T_{eff}) = r_{A}(T^{\diamond}) \leq \sup \{ |\lambda| : \lambda \in \sigma_{A}(T)\}$.
\end{proof}
\begin{remark}
We know that $\|T^{n}\|_{A} = \|(T^{\sharp})^{n}\|_{A}, \text { for all } n \in \mathbb N$ and $T \in B_{A}(H)$. Therefore, $r_{A}(T) = r_{A}(T^{\sharp})$.
\end{remark}

Now, we define an operator $T_{b}$ on $\mathbf{R}(A^{\frac{1}{2}})$ by $T_{b}(A^{\frac{1}{2}} x) = A^{\frac{1}{2}}Tx= (T^{\diamond})^{*}A^{\frac{1}{2}}x, \text{ for all } x \in H$ and $T \in B_{A^{\frac{1}{2}}}(H)$.
\begin{lemma}\label{lemma 0.3.6}
$T_{b}$ is bounded on $\mathbf{R}({A^{\frac{1}{2}}})$ with $T_{b}^{*}= (A^{\frac{1}{2}} T^{*} (A^{\frac{1}{2}})^{\dagger})_{ {R(A^{\frac{1}{2}})}}$.
\end{lemma}
\begin{proof}
$\|T_{b}(A^{\frac{1}{2}} x)\|_{\mathbf{R}({A^{\frac{1}{2}}})} = \|A^{\frac{1}{2}}Tx\|_{\mathbf{R}({A^{\frac{1}{2}}})} = \|PTPx\| \leq \|PT\| \|A^{\frac{1}{2}}x\|_{\mathbf{R}({A^{\frac{1}{2}}})}$. Thus, $T_{b}$ is bounded.

Now, consider an arbitrary element $A^{\frac{1}{2}}z$ in ${\mathbf{R}({A^{\frac{1}{2}}})}$. Then $(T_{b}(A^{\frac{1}{2}}x), A^{\frac{1}{2}}z)= (A^{\frac{1}{2}}Tx, A^{\frac{1}{2}}z)= \langle PTx, Pz\rangle = \langle Px, PT^{*}Pz\rangle = (A^{\frac{1}{2}}x, A^{\frac{1}{2}}T^{*}(A^{\frac{1}{2}})^{\dagger} A^{\frac{1}{2}}z), \text{ for all } x\in H $. Therefore, $T_{b}^{*}= (A^{\frac{1}{2}} T^{*} (A^{\frac{1}{2}})^{\dagger})_{ {R(A^{\frac{1}{2}})}}$.
\end{proof}
\begin{lemma}
Let $T \in B_{A^{\frac{1}{2}}}(H)$. Then $\|T_{b}\|_{\mathbf{R}({A^{\frac{1}{2}}})} = \|(T^{\diamond})^{\diamond}\| = \|T^{\diamond}\|_{A}$.
\end{lemma}
\begin{proof}
 By Proposition \ref{ pro 2.5}, ${\|T\|}_{A} = \|A^{\frac{1}{2}}T(A^{\frac{1}{2}})^{\dagger}\| = \|(A^{\frac{1}{2}})^{\dagger} T^{*}A^{\frac{1}{2}}\| = \|T^{\diamond}\|$. Now, we consider $T^{\diamond}$ instead of $T$ to get $\|(T^{\diamond})^{\diamond}\| = \|T^{\diamond}\|_{A}$. Thus, the equality on the right-hand side is true.

 $\|T_{b}(A^{\frac{1}{2}}x)\|_{\mathbf{R}({A^{\frac{1}{2}}})} = \|PTPx\|= \|(T^{\diamond})^{\diamond} Px\| \leq \|(T^{\diamond})^{\diamond}\| \|A^{\frac{1}{2}}x\|_{\mathbf{R}({A^{\frac{1}{2}}})}, \text{ for all } x\in H$.
 Again, $\|(T^{\diamond})^{\diamond}x\| =\|PTx\| = \|T_{b}(A^{\frac{1}{2}}x)\|_{\mathbf{R}({A^{\frac{1}{2}}})} \leq \|T_{b}\|_{\mathbf{R}({A^{\frac{1}{2}}})} \|x\|, \text{ for all } x\in H$. Therefore, the left-hand equality  $\|T_{b}\|_{\mathbf{R}({A^{\frac{1}{2}}})} = \|(T^{\diamond})^{\diamond}\|$ is also true.
 
\end{proof} 
\begin{theorem}\label{thm 0.3.7}
Let $T \in B_{A^{\frac{1}{2}}}(H)$ and $TA^{\frac{1}{2}} = A^{\frac{1}{2}}T$. Then $\sup \{|\lambda|  : \lambda \in \sigma_{A}(T^{\diamond})\} = r_{A}(T^{\diamond})$.
\end{theorem}
\begin{proof}
Let us consider $|\lambda| > \|T^{\diamond}\|_{A} = \|T_{b}\|_{\mathbf{R}({A^{\frac{1}{2}}})}$. There exists a bounded operator $V_{\lambda} \in B({\mathbf{R}({A^{\frac{1}{2}}})})$ such that 
\begin{align}\label{eq .4}
(\overline{\lambda}- T_{b})V_{\lambda} = V_{\lambda}(\overline{\lambda} - T_{b}) = I_{\mathbf{R}({A^{\frac{1}{2}}})}.
\end{align}
We will show that $\lambda \in \rho_{A}(T^{\diamond})$.
From the equality (\ref{eq .4}), we get
\begin{align}
 V_{\lambda}^{*}({\lambda} - T_{b}^{*}) = ({\lambda}- T_{b}^{*})V_{\lambda}^{*} = I_{\mathbf{R}({A^{\frac{1}{2}}})}.
\end{align}
Considering an arbitrary element $A^{\frac{1}{2}}x\in R(A^{\frac{1}{2}})$, we have
\begin{align}
V_{\lambda}^{*}(\lambda A^{\frac{1}{2}}x - A^{\frac{1}{2}}T^{*}Px) = (\lambda - A^{\frac{1}{2}}T^{*}(A^{\frac{1}{2}})^{\dagger}) V_{\lambda}^{*}A^{\frac{1}{2}}x = A^{\frac{1}{2}}x.
\end{align}
Since, $T$ commutes with $A$. So, $\|PT^{*}\| = \|TP\| = \|PT\| = \|(T^{\diamond})^{\diamond}\| < |\lambda|$ and 
\begin{align*}
V_{\lambda}^{*}(Ax) = \lambda( Ax + \frac{T_{b}^{*}Ax}{\lambda} + \frac{{T_{b}^{*}}^{2}Ax}{\lambda^{2}}+ \dots)= \lambda (Ax + \frac{APT^{*}x}{\lambda} + \frac{A(PT^{*})^{2}x}{\lambda^{2}}+ \dots)
\end{align*}
Now, we can confirm that $V_{\lambda}^{*}(R(A)) \subset R(A)$. It is true that $(A^{\frac{1}{2}})^{\dagger}Z_{A}^{*}V_{\lambda}^{*}W_{A}$ is in $B(H)$ because it is closed in domain $H$. Again, $((A^{\frac{1}{2}})^{\dagger}Z_{A}^{*}V_{\lambda}^{*}W_{A})^{*} = \overline{W_{A}^{*}V_{\lambda} Z_{A}(A^{\frac{1}{2}})^{\dagger}}$. So, $R(((A^{\frac{1}{2}})^{\dagger}Z_{A}^{*}V_{\lambda}^{*}W_{A})^{*} A^{\frac{1}{2}})  \subset R(A^{\frac{1}{2}})$. Thus,  $(A^{\frac{1}{2}})^{\dagger}Z_{A}^{*}V_{\lambda}^{*}W_{A} \in B_{A^{\frac{1}{2}}}(H)$.

Moreover, for all $x\in H$,
\begin{align*}
A (A^{\frac{1}{2}})^{\dagger}Z_{A}^{*}V_{\lambda}^{*}W_{A}(\lambda -T^{\diamond})x= V_{\lambda}^{*}A^{\frac{1}{2}}(\lambda -T^{*}P)A^{\frac{1}{2}}x = Ax
\end{align*}
and 
\begin{align*}
A(\lambda -T^{\diamond})(A^{\frac{1}{2}})^{\dagger}Z_{A}^{*}V_{\lambda}^{*}W_{A}x= (\lambda - A^{\frac{1}{2}}T^{*}(A^{\frac{1}{2}})^{\dagger}) V_{\lambda}^{*}Ax = Ax.
\end{align*}
Hence, $(A^{\frac{1}{2}})^{\dagger}Z_{A}^{*}V_{\lambda}^{*}W_{A}$ is the $A$-inverse of $(\lambda -T^{\diamond})$ in $B_{A^{\frac{1}{2}}}(H)$ and $\lambda \in \rho_{A}(T^{\diamond})$ which implies that $|\lambda| \leq \|T^{\diamond}\|_{A}, \text{ for all }\lambda\in \sigma_{A}(T^{\diamond})$.

 Proposition $5.10$ \cite{MR4362420} gives that 
 \begin{align*}
 |\lambda^{n}| \leq \|(T^{n}) ^{\diamond}\|_{A} = \|(T^{\diamond})^{n}\|_{A}, \text{ for all } \lambda\in \sigma_{A}(T^{\diamond}).
 \end{align*}
 So,  $\sup \{|\lambda|  : \lambda \in \sigma_{A}(T^{\diamond})\} \leq r_{A}(T^{\diamond})$. Therefore, Theorem \ref{thm 5.1} shows that $\sup \{|\lambda|  : \lambda \in \sigma_{A}(T^{\diamond})\} = r_{A}(T^{\diamond})$.
\end{proof}

\begin{theorem}\label{thm  0.3.8}
Let $T \in B_{A^{\frac{1}{2}}}(H)$ and $T A^{\frac{1}{2}} = A^{\frac{1}{2}} T$. Then $r_{A}(T) = \sup \{|\lambda|  : \lambda \in \sigma_{A}(T)\}$.
\end{theorem}
\begin{proof}
Theorem \ref{thm 0.3.7} shows that $r_{A}((T^{\diamond})^{\diamond}) = \sup \{|\lambda|  : \lambda \in \sigma_{A}(PT) = \sigma_{A}((T^{\diamond})^{\diamond})\}$. We know that $r_{A}(T)= r_{A}(PT)= r_{A}((T^{\diamond})^{\diamond})$ and  $\sup \{|\lambda|  : \lambda \in \sigma_{A}(PT)\}= \sup \{|\lambda|  : \lambda \in \sigma_{A}(T)\}$. Therefore, $r_{A}(T) = \sup \{|\lambda|  : \lambda \in \sigma_{A}(T)\}$.
\end{proof}
\begin{corollary}\label{corollary 0.3.9}
Let $T \in B_{A^{\frac{1}{2}}}(H)$ and $T A^{\frac{1}{2}} = A^{\frac{1}{2}} T$. If $|\lambda| > \|T\|_{A}$, then $\lambda \in \rho_{A}(T)$.
\end{corollary}
\begin{proof}
It is obvious that $r_{A}(T) \leq \|T\|_{A}$. By Theorem \ref{thm 0.3.8}, we can show that $\sup \{|\lambda|  : \lambda \in \sigma_{A}(T)\} \leq \|T\|_{A}$. Therefore, for all $\lambda$ with $|\lambda| > \|T\|_{A}$ we get $\lambda \in \rho_{A}(T)$.
\end{proof}
\begin{corollary}
Let $T \in B_{A^{\frac{1}{2}}}(H)$. Assume the condition $T A^{\frac{1}{2}} = A^{\frac{1}{2}} T$ is true . Then $r_{A}(T)= r_{A}(T^{\diamond})$.
\end{corollary}
\begin{proof}
By Lemma (\ref{lemma 0.3.4}) and Theorem (\ref{thm  0.3.8}), we have $r_{A}(T^{\diamond}) \leq r_{A}(T)$. Now we consider $T^{\diamond}$ instead of $T$, we get $r_{A}(T) = r_{A}(PT) = r_{A}((T^{\diamond})^{\diamond}) \leq r_{A}(T^{\diamond})$. Therefore, $r_{A}(T)= r_{A}(T^{\diamond})$.
\end{proof}

Now, we investigate some properties when the condition $T A^{\frac{1}{2}} = A^{\frac{1}{2}} T$ is waned.
\begin{theorem}\label{thm .15}
Let $T \in B_{A^{\frac{1}{2}}}(H)$. Assume both the conditions $\|T\|_{A} < |\lambda|$ and $\|T^{\diamond}\|_{A} < |\lambda|$ hold. Then $\lambda \in \rho_{A}(T^{\diamond})$.
\end{theorem}
\begin{proof}
Let us consider $|\lambda| > \|T^{\diamond}\|_{A} = \|T_{b}\|_{\mathbf{R}({A^{\frac{1}{2}}})}$. There exists a bounded operator $V_{\lambda} \in B({\mathbf{R}({A^{\frac{1}{2}}})})$ such that 
\begin{align}\label{eq .7}
(\overline{\lambda}- T_{b})V_{\lambda} = V_{\lambda}(\overline{\lambda} - T_{b}) = I_{\mathbf{R}({A^{\frac{1}{2}}})}.
\end{align}
We will show that $\lambda \in \rho_{A}(T^{\diamond})$.
From the equality (\ref{eq .7}), we get
\begin{align}
 V_{\lambda}^{*}({\lambda} - T_{b}^{*}) = ({\lambda}- T_{b}^{*})V_{\lambda}^{*} = I_{\mathbf{R}({A^{\frac{1}{2}}})}.
\end{align}
Considering an arbitrary element $A^{\frac{1}{2}}x\in R(A^{\frac{1}{2}})$, we have
\begin{align}
V_{\lambda}^{*}(\lambda A^{\frac{1}{2}}x - A^{\frac{1}{2}}T^{*}Px) = (\lambda - A^{\frac{1}{2}}T^{*}(A^{\frac{1}{2}})^{\dagger}) V_{\lambda}^{*}A^{\frac{1}{2}}x = A^{\frac{1}{2}}x.
\end{align}
We can confirm that $V_{\lambda}^{*}(R(A)) \subset R(A)$ because $\| T^{\diamond}\| = \|T\|_{A} < |\lambda|$ and 
\begin{align*}
V_{\lambda}^{*}(Ax) = \lambda( Ax + \frac{T_{b}^{*}Ax}{\lambda} + \frac{{T_{b}^{*}}^{2}Ax}{\lambda^{2}}+ \dots)= \lambda (Ax + \frac{AT^{\diamond}x}{\lambda} + \frac{A(T^{\diamond})^{2}x}{\lambda^{2}}+ \dots)
\end{align*}
It is true that $(A^{\frac{1}{2}})^{\dagger}Z_{A}^{*}V_{\lambda}^{*}W_{A}$ is in $B(H)$ because it is closed in domain $H$. Again, $((A^{\frac{1}{2}})^{\dagger}Z_{A}^{*}V_{\lambda}^{*}W_{A})^{*} = \overline{W_{A}^{*}V_{\lambda} Z_{A}(A^{\frac{1}{2}})^{\dagger}}$. So, $R(((A^{\frac{1}{2}})^{\dagger}Z_{A}^{*}V_{\lambda}^{*}W_{A})^{*} A^{\frac{1}{2}})  \subset R(A^{\frac{1}{2}})$.\\ 
Thus,  $(A^{\frac{1}{2}})^{\dagger}Z_{A}^{*}V_{\lambda}^{*}W_{A} \in B_{A^{\frac{1}{2}}}(H)$.

Moreover, for all $x\in H$,
\begin{align*}
A (A^{\frac{1}{2}})^{\dagger}Z_{A}^{*}V_{\lambda}^{*}W_{A}(\lambda -T^{\diamond})x= V_{\lambda}^{*}A^{\frac{1}{2}}(\lambda -T^{*}P)A^{\frac{1}{2}}x = Ax
\end{align*}
and 
\begin{align*}
A(\lambda -T^{\diamond})(A^{\frac{1}{2}})^{\dagger}Z_{A}^{*}V_{\lambda}^{*}W_{A}x= (\lambda - A^{\frac{1}{2}}T^{*}(A^{\frac{1}{2}})^{\dagger}) V_{\lambda}^{*}Ax = Ax.
\end{align*}
Therefore, $(A^{\frac{1}{2}})^{\dagger}Z_{A}^{*}V_{\lambda}^{*}W_{A}$ is the $A$-inverse of $(\lambda -T^{\diamond})$ in $B_{A^{\frac{1}{2}}}(H)$ and $\lambda \in \rho_{A}(T^{\diamond})$.

\end{proof}
\begin{corollary}\label{cor .16}
Let $T \in B_{A^{\frac{1}{2}}}(H)$. Then $|\lambda| \leq max\{ \|T\|_{A} , \|T^{\diamond}\|_{A}\}, \text{ for all } \lambda\in \sigma_{A}(T)$. Moreover $\sigma_{A}(T)$ is bounded.
\end{corollary}
\begin{proof}
This statement can be easily proven from Theorem \ref{thm .15}.
\end{proof}
In the next corollary, we present the modified version of Theorem 5.3 \cite{MR4362420}
\begin{corollary}\label{cor .1}
Let $T\in B_{A^{\frac{1}{2}}}(H)$ be invertible in $B_{A^{\frac{1}{2}}}(H)$ with an $A$-inverse $S$ and let $T^{\prime}\in B_{A^{\frac{1}{2}}}(H)$ be such that $\|T^{\prime}S\|_{A} < 1 \text{ and } \|(T^{\prime}S)^{\diamond}\|_{A} < 1 $. Then $T + T^{\prime}$ is $A$-invertible in $B_{A^{\frac{1}{2}}}(H)$.
\end{corollary}
\begin{proof}
The proof follows from Corollary \ref{cor .16} and the analogous proof of  Theorem 5.3 \cite{MR4362420}.
\end{proof}
\begin{theorem}
    Let $T \in B_{A^{\frac{1}{2}}}(H)$. The $A$-resolvent set $\rho_{A}(T)$ is open in $\mathbb{C}$.Moreover there is an analytic function $R: \rho_{A}(T) \mapsto B_{A^{\frac{1}{2}}}(H) $ satisfying $R(\lambda)$ is an $A$-inverse of $\lambda -T$ for all $\lambda \in \rho_{A}(T)$.
\end{theorem}
\begin{proof}
     The theorem can be proven analogously from the Theorem 5.6 \cite{MR4362420} but we have to consider $\min\{{\frac{1}{\|R_{0}\|}, {\frac{1}{\|R_{0}\|_{A}}}}, {\frac{1}{\|R_{0}^{\diamond}\|_{A}}}\}$ instead of $\min\{{\frac{1}{\|R_{0}\|}, {\frac{1}{\|R_{0}\|_{A}}}}\}$ from Theorem 5.6 \cite{MR4362420}.
\end{proof}
\begin{theorem}\label{thm .19}
Let $T \in B_{A^{\frac{1}{2}}}(H)$. Then $\sup\{ |\lambda| : \lambda \in \sigma_{A}(T)\} = \frac{r_{A}(T) + r_{A}(T^{\diamond}) + |r_{A}(T^{\diamond}) - r_{A}(T)|}{2}$. 
\end{theorem}
\begin{proof}
We confirm from Corollary \ref{cor .16} and the first part of Proposition 5.10 \cite{MR4362420} that $\sup\{ |\lambda| : \lambda \in \sigma_{A}(T)\} \leq max\{ \|T^{n}\|_{A}^{\frac{1}{n}}, \|(T^{\diamond})^{n}\|_{A}^{\frac{1}{n}}\}, \text { for all } n\in \mathbb N$.
We know,
\begin{align}\label{equ .17}
max\{ \|T^{n}\|_{A}^{\frac{1}{n}}, \|(T^{\diamond})^{n}\|_{A}^{\frac{1}{n}}\} = \frac{\|T^{n}\|_{A}^{\frac{1}{n}} + \|(T^{\diamond})^{n}\|_{A}^{\frac{1}{n}} +|\|T^{n}\|_{A}^{\frac{1}{n}} - \|(T^{\diamond})^{n}\|_{A}^{\frac{1}{n}}|}{2}.
\end{align}
Now consider limiting both sides of the equation (\ref{equ .17}); we get 
\begin{align*}
\sup\{ |\lambda| : \lambda \in \sigma_{A}(T)\} \leq \frac{r_{A}(T) + r_{A}(T^{\diamond}) + |r_{A}(T^{\diamond}) - r_{A}(T)|}{2}.
\end{align*}
Again, $\sup\{ |\lambda| : \lambda \in \sigma_{A}(T)\} \geq \frac{r_{A}(T) + r_{A}(T^{\diamond}) + |r_{A}(T^{\diamond}) - r_{A}(T)|}{2}$ follows from Theorem \ref{thm 5.1} and Lemma \ref{lemma 0.3.4}. Therefore, 
\begin{align*}
\sup\{ |\lambda| : \lambda \in \sigma_{A}(T)\} = \frac{r_{A}(T) + r_{A}(T^{\diamond}) + |r_{A}(T^{\diamond}) - r_{A}(T)|}{2}.
\end{align*}

\end{proof}
\begin{corollary}
Let $T \in B_{A^{\frac{1}{2}}}(H)$. Then either $r_{A}(T) = \sup\{|\lambda| : \lambda \in \sigma_{A}(T)\}$ or  $r_{A}(T^{\diamond}) = \sup\{|\lambda| : \lambda \in \sigma_{A}(T)\}$.
\end{corollary}
\begin{proof}
If $r_{A}(T) \geq r_{A}(T^{\diamond})$, then Theorem \ref{thm .19} and Theorem \ref{thm 5.1} say that $r_{A}(T) = \sup\{|\lambda| : \lambda \in \sigma_{A}(T)\}$. Similarly, when $r_{A}(T^{\diamond}) \geq r_{A}(T)$, then Theorem \ref{thm .19} and Lemma \ref{lemma 0.3.4} confirm that $r_{A}(T^{\diamond}) = \sup\{|\lambda| : \lambda \in \sigma_{A}(T)\}$.
\end{proof}
$\bold{ Note  }$: If it will be possible to find an example $ r_{A}(T_{1}^{\diamond}) > r_{A}(T_{1}) \text{ for some } T_{1} \in B_{A^{\frac{1}{2}}}(H)$ then the open problem ($\sup\{ |\lambda|: \lambda \in \sigma_{A}(T)\} = r_{A}(T) , \text{ when }T \in B_{A^{\frac{1}{2}}}(H)$) will have a negative answer. Otherwise,  $r_{A}(T) \geq r_{A}(T^{\diamond})  \text{ and } r_{A}(T) \leq r_{A}(T^{\diamond})  \text { for all } T \in B_{A^{\frac{1}{2}}}(H)$ say that $r_{A}(T^{\diamond}) = r_{A}(T) = \sup\{|\lambda| : \lambda \in \sigma_{A}(T)\}, \text{ for all } T \in B_{A^{\frac{1}{2}}}(H)$. More precisely, it is enough to show that   $r((T_{1}^{\diamond})^{\diamond}) > r(T_{1}^{\diamond}) , \text { for some } T_{1} \in B_{A^{\frac{1}{2}}}(H)$ in order to disprove that open problem. 

Now, we justify our proven results by following examples.

\begin{example}
Let $A$ be a diagonal operator in ${\ell_{\mathbb Z}^{2}}(\mathbb C)$ given by $A(e_{-n})= {\frac{1}{n}} e_{-n}$, $A(e_{n})= {\frac{1}{n^{2}}} e_{n}, \text{ for all } n\in \mathbb N \setminus \{1\}$ and $A(e_{-1}) = A(e_{1}) = A(e_{0}) = 0$. Define $T: {\ell_{\mathbb Z}^{2}}(\mathbb C) \mapsto {\ell_{\mathbb Z}^{2}}(\mathbb C)$ by $T(e_{n}) = \frac{1}{\sqrt{n}} e_{-n}, T(e_{-n}) = e_{-n} \text{ and } T(e_{0})= 0$. It is easy to show that $T \text { is in } B_{A^{\frac{1}{2}}}(H)$ but not in $B_{A}(H)$. Also,  $r_{A}(T) = 1$ because $\|T^{k}\|_{A} = 1, \text { for all } k\in \mathbb N$. Moreover, $T^{\diamond}e_{n} = T^{\diamond}e_{-1} = T^{\diamond}(e_{0})= 0, \text{ for all } n\in\mathbb N$ and $T^{\diamond}(e_{-n}) = e_{-n} + e_{n}, \text { for all } n\geq 2$. Then $\|(T^{\diamond})^{k}\|_{A} = \sqrt{\frac{3}{2}}, \text{ for all } k \in \mathbb N$. So, $r_{A}(T) = r_{A}(T^{\diamond}) = 1$. We can also find that $\rho_{A}(T) = \mathbb C \setminus\{1\}$. Thus, $\sigma_{A}(T) = \{1\}$. Therefore, $r_{A}(T) = r_{A}(T^{\diamond}) = \sup\{|\lambda| : \lambda \in \sigma_{A}(T)\} = 1$.
\end{example}
\begin{example}
Let $A$ be a diagonal operator in ${\ell_{\mathbb Z}^{2}}(\mathbb C)$ given by $A(e_{-n})= {\frac{1}{n}} e_{-n}$, $A(e_{n})= {\frac{1}{n^{2}}} e_{n}, \text{ for all } n\in \mathbb N \setminus \{1\}$ and $A(e_{-1}) = A(e_{1}) = A(e_{0}) = 0$. Define $T: {\ell_{\mathbb Z}^{2}}(\mathbb C) \mapsto {\ell_{\mathbb Z}^{2}}(\mathbb C)$ by $T(e_{n}) = \frac{1}{\sqrt{n}} e_{-n},  T(e_{-n}) = e_{n} \text{ and } T(e_{0})= 0$. It is easy to show that $T \text { is in } B_{A^{\frac{1}{2}}}(H)$ but not in $B_{A}(H)$. Moreover, $T^{\diamond}(e_{-n}) = e_{n}, T^{\diamond}e_{n} = {\frac{1}{\sqrt n}} e_{-n}, \text{ for all } n \geq 2 \text{ and } T^{\diamond}(e_{0}) = T^{\diamond}(e_{1}) = T^{\diamond}(e_{-1}) =0$. We can show that $r_{A}(T) = r_{A}(T^{\diamond}) = \frac{1}{\sqrt[4]2}$. Again, $\sigma_{A}(T) = \{\lambda \in \mathbb C : {\lambda}^{2} = \frac{1}{\sqrt n}, \text{ for all } n \in \mathbb N\}$. Therefore, $r_{A}(T) = r_{A}(T^{\diamond}) = \sup\{|\lambda| : \lambda \in \sigma_{A}(T)\} = \frac{1}{\sqrt[4]2}$.

\end{example}
 
\begin{example}
Let $T : \ell^{2} \to \ell^{2}$ be the left shift operator defined by  $T(x_{1}, x_{2}, x_{3}, \dots, x_{n}, \dots) = (x_{2}, x_{3},\dots, x_{n}, \dots)$  and $A : \ell^{2} \to \ell^{2}$ be defined by $$A(x_{1}, x_{2}, x_{3}, \dots, x_{n}, \dots) = \Big(x_{1}, \frac{x_{2}}{2^{2}}, \frac{x_{3}}{2^{4}},\dots, \frac{x_{n}}{2^{2(n-1)}}, \dots\Big).$$ Then, $A$ is a positive operator. Moreover, $T\in B_{A}(H)$. Also,  $r_{A}(T) = 2$ because $\|T^{n}\|_{A} = 2^{n}, \text{ for all } n\in \mathbb N$. It is easy to show that $\lambda \in \rho_{A}(T)$ if and only if  $2 < |\lambda|$. Thus, $|\lambda| \leq 2$ if and only if  $\lambda \in \sigma_{A}(T)$. Therefore, $r_{A}(T) = \sup \{|\lambda| : \lambda \in \sigma_{A}(T)\} = 2$.
\end{example}
\begin{theorem}
 Let $T \in B_{A^{\frac{1}{2}}}(H)$ be $A$-invertible in $B_{A^{\frac{1}{2}}}(H)$ with an $A$-inverse $S$ in $B_{A^{\frac{1}{2}}}(H)$. Then $r_{A}(T)r_{A}(S) \geq 1$.
 \end{theorem}
 \begin{proof}
We know that $T^{*}A^{\frac{1}{2}} = A^{\frac{1}{2}} T^{\diamond}$. It can be observed that $N(A^{\frac{1}{2}} T^{\diamond}) \neq H$. Otherwise, $T^{*}A^{\frac{1}{2}} = 0  \text{ if and only if } AT = 0$. So, $A= ATS = 0$ is a contradiction. There exists an element $x_{0} \in H \setminus N(A^{\frac{1}{2}}T^{\diamond})$ such that $0 \neq \|T^{\diamond}x_{0}\|_{A} = \|STT^{\diamond} x_{0}\|_{A} \leq \|S\|_{A} \|T\|_{A} \|T^{\diamond}x_{0}\|_{A} $. Thus, $\|S\|_{A} \|T\|_{A} \geq 1$. By using induction we get $S^{n}$ is the $A$-inverse in $B_{A^{\frac{1}{2}}}(H)$ of $T^{n}$, $\text{ for all } n\in \mathbb N$. So, $\|S^{n}\|_{A}^{\frac{1}{n}} \|T^{n}\|_{A}^{\frac{1}{n}} \geq 1, \text{ for all }n \in \mathbb N$. Therefore, $r_{A}(S)r_{A}(T) \geq 1$.
 \end{proof}
 \begin{theorem}
 Let $T \in B_{A^{\frac{1}{2}}}(H)$. $T$ is an $A$-invertible operator in $B_{A^{\frac{1}{2}}}(H)$ if and only if  there exist two operator $S_{1}$ and $S_{2}$ in $B_{A^{\frac{1}{2}}}(H)$ such that $ATS_{1}= AS_{2}T = A$.
 \end{theorem}
 \begin{proof}
The first part is obvious from the definition of $A$-invertible operator in $B_{A^{\frac{1}{2}}}(H)$.

Conversely, we claim that $N(AT) = N(A)$. Consider an element $x \in N(AT)$, we get $Ax = AS_{2}PTx= 0$. So $N(AT) \subset N(A)$. Similarly, taking $z\in N(A) =N(P)$, we get $ATz = APTPz= 0$. Thus, $N(AT) = N(A)$. By Proposition(3.7) \cite{MR4362420}, we confirm that $S_{1}$ is an $A$-inverse of $T$.
\end{proof}
\section{Characterization of $A$-Harte spectral radius for commuting  $A^{\frac{1}{2}}$-adjoint operators}
 Let $\bold{T} = (T_{1},\dots, T_{d}) \in {B(H)}^{d}$ be an $d$-tuple of pairwise commuting operators on $H$. The symbol $\sigma_{h}(\bold{T})$ will stand for the Harte spectrum of $T$, i.e. $(\lambda_{1},\dots, \lambda_{d}) \notin \sigma_{h}(\bold{T})$ if there exist operators $U_{1},\dots,U_{d}$ and $V_{1},\dots,V_{d}$ in $B(H)$ such that $\sum_{j=1}^{d} U_{j}(\lambda_{j} - T_{j}) = I$ and $\sum_{j=1}^{d}(\lambda_{j} -T_{j})V_{j}= I$. Here, the Harte spectral radius of $\bold{T}$ is defined as $r_{h}(\bold{T}) = \sup \{\|\lambda\|_{2}: \lambda = (\lambda_{1},\dots, \lambda_{d}) \in \sigma_{h}(\bold{T})\}$.
 \begin{theorem}\cite{MR1202017}
 Let $\bold{T} = (T_{1},\dots,T_{d})$ be a mutually commuting $d$-tuple of Hilbert space operators. Then
 \begin{align*}
 r_{h}(\bold{T}) = \lim_{n \to \infty} \|\sum_{|s|=n, s\in \mathbb{Z}_{+}^{d}} \frac{n!}{s!} ({\bold {T}}^{*})^{s}\bold {T}^{s}\|^{\frac{1}{2n}}\\
 = {\inf_{n \in \mathbb{N}}}~ \|\sum_{|s|=n, s\in \mathbb{Z}_{+}^{d}} \frac{n!}{s!} (\bold {T}^{*})^{s}\bold{T}^{s}\|^{\frac{1}{2n}}
 \end{align*}
 
 \end{theorem}
 \begin{definition}
 Let $\bold{T}=(T_{1},\dots,T_{d}) \in (B_{A^{\frac{1}{2}}}(H))^{d}$ be a mutually commuting $d$-tuple of $B_{A^{\frac{1}{2}}}(H)$ operators. Then $\lambda =(\lambda_{1},\dots,\lambda_{2})$ is called a point of A-Harte resolvent set of $\bold{T}$ if there exist $(U_{1},\dots,U_{d}) \in (B_{A^{\frac{1}{2}}}(H))^{d}$ and $(V_{1},\dots,V_{d}) \in (B_{A^{\frac{1}{2}}}(H))^{d}$ such that  $P\sum_{j=1}^{d} U_{j}(\lambda_{j} - T_{j}) = P$ and $P\sum_{j=1}^{d}(\lambda_{j} -T_{j})V_{j}= P$. 

 Here, $\rho_{A_{h}}(\bold{T})$ denotes $A$-Harte resolvent set of $\bold{T}$ and $\sigma_{A_{h}}(\bold{T}) = \mathbb{C}^{d} \setminus \rho_{A_{h}}(\bold{T})$ is called $A$-Harte spectrum of $\bold{T}$ . 
 \end{definition}
 \begin{definition}
 Let $\bold{T}=(T_{1},\dots,T_{d}) \in (B_{A^{\frac{1}{2}}}(H))^{d}$ be a mutually commuting $d$-tuple. $A$-Harte spectral radius of $\bold{ T }$ is defined as 
 \begin{align}
 r_{A_{h}}(\bold{T}) = \lim_{n \to \infty} \|\sum_{|s|=n, s\in \mathbb{Z}_{+}^{d}} \frac{n!}{s!} (\bold{T}^{\diamond})^{s} \bold{T}^{s}\|_{A}^{\frac{1}{2n}}\\
 = {\inf_{n \in \mathbb{N}}}~ \|\sum_{|s|=n, s\in \mathbb{Z}_{+}^{d}} \frac{n!}{s!} (\bold{T}^{\diamond})^{s} \bold{T}^{s}\|_{A}^{\frac{1}{2n}}
 \end{align}
 
 \end{definition}
 Now, it is natural to raise the question whether $r_{A_{h}}(\bold{T}) = \sup \{\|\lambda\|_{2}: \lambda \in \sigma_{A_{h}}(\bold{T})\}$ is true or not. We prove that one side inequality  $r_{A_{h}}(\bold{T}) \leq \sup \{\|\lambda\|_{2}: \lambda \in \sigma_{A_{h}}(\bold{T})\}$  when each $T_{i} ( i= 1,2,\dots,d)$ commutes with $A$. An interesting question remains open to prove whether the opposite side inequality $r_{A_{h}}(\bold{T}) \geq \sup \{\|\lambda\|_{2}: \lambda \in \sigma_{A_{h}}(\bold{T})\}$ is true or not when each $T_{i}$ commutes with $A$.

 We define a new operator $\bold{\hat{T}} = ({\hat{T_{1}},\dots,\hat{T_{d}}})$ in $(B(\bold{R(A^{\frac{1}{2}}})))^{d}$ by $\hat{T}_{i}(A^{\frac{1}{2}} x) = A^{\frac{1}{2}} T_{i}x, \text{ where } ~i= 1,2,\dots,d$ and $\bold{T}=(T_{1},\dots,T_{d}) \in (B_{A^{\frac{1}{2}}}(H))^{d}$. Lemma \ref{lemma 0.3.6} says that $(\hat{T}_{i})^{*}= (A^{\frac{1}{2}} T_{i}^{*} (A^{\frac{1}{2}})^{\dagger})_{ {R(A^{\frac{1}{2}})}}, \text{ for all } i \in \{1,2,\dots,d\}$.
 \begin{theorem}
 Let $\bold{T}=(T_{1},\dots,T_{d}) \in (B_{A^{\frac{1}{2}}}(H))^{d}$ be a mutually commuting $d$-tuple and each $T_{i}, i= 1,2,\dots,d$ commutes with $A$. Then $r_{A_{h}}(\bold{T}) = r_{h}(\bold{\hat{T}})$.
 \end{theorem}
 \begin{proof}
 It is easy to prove that $(\hat{T_{1}},\dots, \hat{T_{d}})$ is a mutually commuting $d$-tuple. For each $S \in B_{A^{\frac{1}{2}}}(H)$,  $(W_{A}^{*})^{\dagger} S^{*} W_{A}^{*}$ is bounded. When two arbitrary operators $V, W \in B_{A^{\frac{1}{2}}}(H)$, then
 \begin{align*}
 & \alpha(V + W)\\
 &= \overline{W_{A}(V + W) W_{A}^{\sharp}}\\
 &=((W_{A}^{*})^{\dagger} (V^{*} + W^{*}) W_{A}^{*})^{*}\\
 &= \alpha(V) + \alpha(W).
\end{align*}
We know $\alpha$ is a homomorphism. Then,
\begin{align}\label{eqn 13}
\begin{split}
&\| \sum_{|s|=n, s \in \mathbb{Z}_{+}^{d}} \frac{n!}{s!} (\bold{T}^{\diamond})^{s} \bold{T}^{s}\|_{A}^{\frac{1}{2n}}\\
&=\| \sum_{|s|=n, s \in \mathbb{Z}_{+}^{d}} \frac{n!}{s!}\alpha(\bold{T}^{\diamond})^{s} \alpha(\bold{T}^{s})\|_{\bold{R(A^{\frac{1}{2}}})}^{\frac{1}{2n}}\\
&=\| \sum_{|s=(s_{1},\dots,s_{d})|=n, s \in \mathbb{Z}_{+}^{d}} \frac{n!}{s!} \overline{W_{A}(T_{1}^{\diamond})^{s_{1}}\dots (T_{d}^{\diamond})^{s_{d}} T_{1}^{s_{1}}\dots T_{d}^{s_{d}}W_{A}^{\sharp}}\|_{\bold{R(A^{\frac{1}{2}}})}^{\frac{1}{2n}}\\
&=\| \sum_{|s=(s_{1},\dots,s_{d})|=n, s \in \mathbb{Z}_{+}^{d}} \frac{n!}{s!} {W_{A}(T_{1}^{\diamond})^{s_{1}}\dots (T_{d}^{\diamond})^{s_{d}} T_{1}^{s_{1}}\dots T_{d}^{s_{d}}W_{A}^{\sharp}}\|_{\bold{R(A^{\frac{1}{2}}})}^{\frac{1}{2n}}\\
\end{split}
\end{align}
Now, 
\begin{equation}\label{eqn 14}
{W_{A}(T_{1}^{\diamond})^{s_{1}}\dots (T_{d}^{\diamond})^{s_{d}} T_{1}^{s_{1}}\dots T_{d}^{s_{d}}W_{A}^{\sharp}}(Ax)\\
=\hat{(T_{1}}^{*})^{s_{1}}\dots (\hat{T_{d}}^{*})^{s_{d}} \hat{T_{1}}^{s_{1}}\dots \hat{T_{d}}^{s_{d}}(Ax)
\end{equation}
$R(A)$ is dense in $\bold{R(A^{\frac{1}{2}}})$. From Equations (\ref{eqn 13}) and (\ref{eqn 14}), we get
\begin{align*}
&\| \sum_{|s|=n, s \in \mathbb{Z}_{+}^{d}} \frac{n!}{s!} (\bold{T}^{\diamond})^{s} \bold{T}^{s}\|_{A}^{\frac{1}{2n}}\\
&=\| \sum_{|s|=n, s \in \mathbb{Z}_{+}^{d}} \frac{n!}{s!} (\bold{\hat{T}}^{*})^{s} \bold{\hat{T}}^{s}\|_{\bold{R(A^{\frac{1}{2}}})}^{\frac{1}{2n}}\\
\end{align*}
Therefore,
$r_{A_{h}}(\bold{T}) = r_{h}(\bold{\hat{T}})$.
\end{proof}
\begin{theorem}
Let $\bold{T}=(T_{1},\dots,T_{d}) \in (B_{A^{\frac{1}{2}}}(H))^{d}$ be a mutually commuting $d$-tuple and each $T_{i}, i= 1,2,\dots,d$ commutes with $A$. Then, $r_{A_{h}}(\bold{T}) \leq \sup\{\|\mu\|_{2}: \mu= (\mu_{1},\dots,\mu_{d}) \in \sigma_{A_{h}}(\bold{T})\}$.
\end{theorem}
\begin{proof}
Consider an arbitrary element $\lambda = (\lambda_{1},\dots,\lambda_{d})\in \rho_{A_{h}}(\bold{T})$, then there exist $(S_{1},\dots,S_{d}) \in (B_{A^{\frac{1}{2}}}(H))^{d}$ and $(V_{1},\dots,V_{d})\in (B_{A^{\frac{1}{2}}}(H))^{d}$ such that
\begin{align*}
P\sum_{i=1}^{d}(\lambda_{i}-T_{i})S_{i}= P  \text{ and }
P\sum_{i=1}^{d}V_{i}(\lambda_{i}- T_{i}) =P.
\end{align*}
So,
\begin{equation}
\sum_{i=1}^{d}(\lambda_{i}- \hat{T_{i}}) Z_{A}S_{i}Z_{A}^{*} = I_{\bold{R(A^{\frac{1}{2}}})}
\end{equation}
and,
\begin{equation}
\sum_{i=1}^{d}(Z_{A}V_{i}Z_{A}^{*})(\lambda_{i} - \hat{T_{i}}) = I_{\bold{R(A^{\frac{1}{2}}})}
\end{equation}
where $Z_{A}S_{i}Z_{A}^{*} \text{ and } Z_{A}V_{i}Z_{A}^{*}$ both are bounded in $H$, for all $i \in \{1,2,\dots,d\}$.
Thus, $\lambda \in \rho_{h}(\bold{\hat{T}})$. So, $\sigma_{h}(\bold{\hat{T}}) \subset \sigma_{A_{h}}(\bold{T})$. 
Therefore, $r_{A_{h}}(\bold{T}) = r_{h}(\bold{\hat{T}}) = \sup\{\|\mu\|_{2}: \mu \in \sigma_{h}(\bold{\hat{T}})\} \leq \sup\{\|\mu\|_{2}: \mu \in \sigma_{A_{h}}(\bold{T}) \}$.  
\end{proof}

\begin{center}
	\textbf{Acknowledgements}
\end{center}

\noindent The present work of the second author was partially supported by Science and Engineering Research Board (SERB), Department of Science and Technology, Government of India (Reference Number: MTR/2023/000471) under the scheme ``Mathematical Research Impact Centric Support (MATRICS)''.


\begin{thebibliography}{10}
	     
      \bibitem{MR2252709}
	Ali Mostafazadeh. \newblock { Metric operator in pseudo-Hermitian quantum mechanics and the imaginary cubic potential}.
	\newblock {\em J. Phys. A,} 39(32):10171--10188, 2006.
 

	   \bibitem{MR1144600}
	F. G. Scholtz, H. B. Geyer, and F. J. W. Hahne. \newblock { Quasi-Hermitian operators in quantum mechanics and the variational principle}.
	\newblock {\em Ann. Physics,} 213(1):74--101, 1992.
      
      \bibitem{MR3834203}
	Hamadi Baklouti, Kais Feki, and Ould Ahmed
              Mahmoud Sid Ahmed.
	\newblock {Joint numerical ranges of operators in semi-Hilbertian spaces}.
	\newblock {\em Linear Algebra Appl.,} 555:266--284, 2018.

      \bibitem{MR4362420}
	 Hamadi Baklouti and Sirine Namouri.\newblock { Spectral analysis of bounded operators on semi-Hilbertian spaces}.
	\newblock {\em Banach J. Math. Anal.,} 16(1):Paper No. 12, 17, 2022.
 
       \bibitem{hb}
	Hamadi Baklouti and Mohamed Mabrouk.\newblock{A note on the A-spectrum of A-bounded operators}.
	\newblock {\em Operators and Matrices,} 17(3):599--611, 2023.

     \bibitem{MR4133606}
	 Kais Feki.
	\newblock {Spectral radius of semi-Hilbertian space operators and its applications}.
	\newblock {\em Ann. Funct. Anal.,} 11(4):929--946, 2020.
 
       \bibitem {MR4402703}
    Kais Feki.
      \newblock {Some {$A$}-spectral radius inequalities for {$A$}-bounded
              {H}ilbert space operators},
       \newblock {Banach J. Math. Anal.}, 16(2):Paper No. 31, 22, 2022.
 
 
		\bibitem{MR2388631}
	 M. Laura Arias, Gustavo Corach, and M. Celeste Gonzalez.
	\newblock {Partial isometries in semi-Hilbertian spaces}.
	\newblock {\em Linear Algebra Appl.,} 428(7):1460--1475, 2008.
 
		\bibitem{MR2442900}
	 M. Laura Arias, Gustavo Corach, and M. Celeste Gonzalez.
	\newblock{Metric properties of projections in semi-Hilbertian spaces}.
	\newblock {\em Integral Equations Operator Theory,} 62(1):11--28, 2008.
 
      \bibitem{MR2590353}
	M. Laura Arias, Gustavo Corach, and M. Celeste Gonzalez.
	\newblock {Lifting properties in operator ranges}.
	\newblock {\em Acta Sci. Math. (Szeged),} 75(3-4):635--653, 2009.
 
    
 
     
 
      \bibitem{MR0203464}
	 R. G. Douglas.
	\newblock {On majorization, factorization, and range inclusion of operators on Hilbert space}.
	\newblock {\em Proc. Amer. Math. Soc.,} 17:413--415, 1966.
 
      \bibitem {MR1202017}
    V. M\"{u}ller, and A. So\l tysiak ,
     \newblock {Spectral radius formula for commuting {H}ilbert space
              operators},
   \newblock {\em Studia Mathematica,} 3(103):329--333, 1992.
      
  

 
  \end{thebibliography}
\end{document}